\documentclass[preprint,12pt]{article}

\usepackage{amssymb}
\usepackage{amscd,amsmath,amsfonts,graphicx}
\usepackage{amsthm}

\newcommand{\R}{\mathbb{R}}

\newcommand{\e}{\rm e}

\date{\today}
\begin{document}

\begin{center}
 {\bf NEW CHARACTERIZATIONS OF THE GAMMA DISTRIBUTION VIA INDEPENDENCE
OF TWO STATISTICS BY USING\\ ANOSOV'S THEOREM }%\footnote{Received by the editors 22 July 2022 \quad Revised version ... }}
\end{center}

\vspace{0.2cm}
\begin{center}
GWO DONG LIN\,\footnote{Institute of Statistical Science, Academia Sinica, Nankang, Taipei 11529, Taiwan (ROC); e-mail:
gdlin@stat.sinica.edu.tw}, \quad JORDAN M. STOYANOV\,\footnote{Institute of Mathematics and Informatics, Bulgarian Academy of Sciences, 1113 Sofia, Bulgaria; e-mail: stoyanovj@gmail.com}
\end{center}

\vspace{0.1cm}
{\small
 \vspace{1cm}\noindent
{\bf Abstract.} Available in the literature are properties which
characterize the gamma distribution via independence of two
appropriately chosen statistics. Well-known is the classical result
when one of the statistics is the sample mean and  the other one the
sample coefficient of variation. In this paper, we elaborate on a
version of Anosov's theorem which allows to establish a general
result, Theorem 1, and a series of seven corollaries providing new
characterization results for gamma distributions. We keep the sample
mean as one of involved statistics, while now the second one  can be
taken from a quite large class of {\it homogeneous feasible definite
statistics}.
 It is relevant to mention that there is an interesting parallel between
 the new characterization results for gamma distributions  and recent characterization results for the normal distribution.

\vspace{0.1cm}
{\bf Key words}:   characterization of distributions, order statistics, gamma distribution, sample mean,
sample coefficient of variation, Anosov's theorem, feasible definite statistics,  normal distribution.

\vspace{0.1cm} {\bf Mathematics Subject Classification 2020:} 60E05,
\ 62E10, \ 62G30

\bigskip
{\bf 1. Introduction.} As usual we write $X \sim \gamma (a, b)$ and say that  $X$ has a {\it gamma distribution} with
parameters $a > 0$ and  $b > 0$ if $X$ is absolutely continuous and has density
\[
f(x)=\frac{1}{\Gamma(a)\,b^a}\,x^{a-1}{\e}^{-x/b},\  \ x>0; \quad f(x) =0, \ x \leq 0.
\]

This distribution is well studied and widely used in both theory and applications.  Our main interests are
properties which uniquely characterize this distribution.

We start with a general situation. Suppose that $X \sim F$ is a
non-degenerate and positive random variable, so the distribution
function $F$ (we use just `a distribution') has a support contained
in ${\mathbb R}_+ :=(0,\infty).$ Let $X_1, X_2, \ldots, X_n$ be a
random (independent) sample of size $n\ge 2$ from
 $F.$ Consider the {\it sample mean}
$\overline{X}_n:=\frac1n\sum_{i=1}^nX_i$, the {\it sample variance}
$S_n^2:=\frac{1}{n-1}\sum_{i=1}^n(X_i-\overline{X}_n)^2,$ and the {\it sample deviation} $S_n = {\sqrt {S_n^2}}$.
Well-known is the following characterization property.

\smallskip
{\bf Classical result:} {\it Under the above conditions and notations, the sample mean $\overline{X}_n$ and the sample
coefficient of variation $S_n/\overline{X}_n$ are independent {\bf \emph{if and only if}}
 $F$ is a gamma distribution.}

 Notice that for this general result there are no a priori moment or smoothness conditions.
We mention briefly that the proof of the sufficiency part is based
on the fact that if $F$ is a gamma distribution, then the sample
mean $\overline{X}_n$ is independent of the random vector
$(X_1/\overline{X}_n,X_2/\overline{X}_n,\ldots, X_n/\overline{X}_n)$
(see [10], [9]). Hence $\overline{X}_n$
is independent of $S_n^2/\overline{X}_n^2,$ because the latter is a
function of the random vector. This in turn implies that
$\overline{X}_n$ and $S_n/\overline{X}_n$ are independent. For the
proof of  the necessity part (see [6], p.\,753).

The paper is organized as follows. In Section 2, we recall the
definition of feasible definite statistics as given and exploited in [3].
In Section 3 we formulate the main results,
Theorem 1 and a series of seven specific and interesting
corollaries. The needed tools and lemmas are given in Section 4. In
Section 5, we prove the main theorem. Finally, in Section 6 we make
some comments and outline two conjectures.

In what follows, we use the symbol \ $\perp\!\!\!\perp$ \ to indicate that some random quantities are `independent'.
The abbreviation `iff' is used for `if and only if', while `a.s.' stands for `almost surely'.

\bigskip
{\bf 2. Feasible definite statistics.}
First we introduce a {\bf \emph{base function}}, $U(\lambda), \lambda \in {\bf A}$, where ${\bf A} \subset \R^n$ is the ordered set:
\begin{eqnarray}
{\bf A}=\Big\{\lambda:=(\lambda_1, \lambda_2,  \ldots, \lambda_n): \ \lambda_1\le \lambda_2\le\cdots\le \lambda_n, \
\Sigma_{i=1}^n\lambda_i=0\Big\}.
\end{eqnarray}

We say that $U$ is  a {\bf \emph{feasible definite function}} of
positive degree $p$ of homogeneity (briefly, of h-degree $p$),
if it satisfies  the following conditions:\\
 (i) \ $U$ is nonnegative and continuous on ${\bf A}$ \ (nonnegativity and continuity);\\
(ii) \ $U(\lambda)=0$ \ iff \ $\lambda=(0, 0, \ldots,0)$ \ (definiteness);\\
(iii) \ $U(s\lambda)=s^p\,U(\lambda)$ for all  $s>0$ and all
$\lambda=(\lambda_1, \lambda_2, \ldots, \lambda_n)\in {\bf A}$ \
(homogeneity

\hspace{0.3cm} of degree $p$).

\smallskip
Suppose that $U(\lambda), \lambda \in {\bf A},$ is a  feasible definite function of h-degree $p$ and let
$X_{(1)}\le X_{(2)}\le \cdots\le X_{(n)}$
be the order statistics of the random sample $X_1, X_2,\ldots, X_n$ of size $n\ge 2$ from $F.$ We define
\begin{eqnarray}
Z_n=U(X_{(1)}-\overline{X}_n,X_{(2)}-\overline{X}_n,\ldots,X_{(n)}-\overline{X}_n),
\end{eqnarray}
and say that $Z_n$ is a {\bf \emph{feasible definite statistic}} on
${\bf A}$ of positive degree $p$ of homogeneity, or, simply, that
$Z_n$ is a feasible definite statistic of h-degree $p.$

It is easy to see that if  $Z_n$ is a feasible definite statistic of
h-degree $p,$ then $(Z_n)^{1/p}$ is a feasible definite statistic of
h-degree $1.$ Notice that the sample mean $\overline{X}_n$ is
however not a feasible definite statistic because it is not of the
form (2).

\bigskip
{\bf 3. Main results: Theorem 1 and Corollaries.} We deal with a random variable $X\sim F$
having {\it positive continuous density} $f$ on $(0,\infty).$ Let
$X_{(1)}\le X_{(2)}\le \cdots\le X_{(n)}$ be the order statistics of
a random sample $X_1, X_2,\ldots, X_n$ of size $n\ge 3$ from
 $F$ and $\overline{X}_n$ its sample mean.
We formulate now our main result.

\medskip
{\bf Theorem 1.}
{\it Suppose that $U(\lambda), \lambda \in {\bf A},$ is a feasible definite function of h-degree 1. Define
 $Z_n=U(X_{(1)}-\overline{X}_n,X_{(2)}-\overline{X}_n,\ldots,X_{(n)}-\overline{X}_n),$
 a feasible definite statistic of h-degree $1$. Then the sample mean $\overline{X}_n$ and the
quotient $V_n:=Z_n/\overline{X}_n$ are independent iff $F$ is a gamma distribution. Symbolically:
\[
\overline{X}_n \ \perp\!\!\!\perp \ V_n:=Z_n/\,\overline{X}_n \quad \Longleftrightarrow \quad F \mbox{ is a gamma distribution}.
\]
}
\indent{\it Remark} 1.  Notice that in the classical result given in
the Introduction, involved are two statistics, the sample mean
$\overline{X}_n$ and the sample coefficient of variation
$S_n/\overline{X}_n.$ In Theorem 1, while keeping one of the
statistics,
 $\overline{X}_n,$ we enjoy the freedom of choosing
the second statistic from a large class. This fact will be well illustrated
 by formulating seven specific and diverse corollaries followed by useful
brief comments. For each corollary it will be clear which base
function $U$ is used.

{\it Remark} 2. Note also that if  instead
of ${\R}_+$ we assume that $X\sim F$ has a {positive continuous
density} $f$ on the {\it whole real line} $\R=(-\infty,\infty),$
then the independence of the sample mean $\overline{X}_n$ and the
homogeneous feasible definite statistic $Z_n$ (rather than $V_n$),
defined in (2) or in Theorem 1, will characterize $F$ as being a
normal distribution. Therefore, our results for gamma distributions
can be considered in a sense as counterparts of the normal
characterizations derived recently in [3].

The statements in Corollaries 1 and 5 are known. They are obtained
by a different approach in [7] %Hwang and Hu (2000)
and are listed here
for completeness of the general picture and for comparison with the
other corollaries which are new. In fact, the known Corollaries 1
and 5 are special cases of our Corollaries 2 and 4, as explained
below.

\smallskip
{\bf Corollary 1.} {\it Let $R_{n}:=X_{(n)}-X_{(1)}$ be
the range of the random sample. Then the sample mean
$\overline{X}_n$ and  the quotient  $R_n/\overline{X}_{n}$ are
independent iff \ $F$ is a gamma distribution, i.e.,
\[
\overline{X}_n \ \perp\!\!\!\perp \ R_n/\,\overline{X}_n \quad \Longleftrightarrow \quad F \mbox{ is a gamma distribution}.
\]
} \indent Indeed, note that
$R_{n}=X_{(n)}-X_{(1)}=(X_{(n)}-\overline{X}_n)-(X_{(1)}-\overline{X}_n)$
is a feasible definite statistic of h-degree 1, so the conditions in
Theorem 1 are satisfied.

\smallskip
{\bf Corollary 2.}
{\it  Let $a_1, a_2, \ldots, a_n$ be arbitrary real numbers, not all coinciding, and such that $a_1 \le a_2 \le
\cdots\le a_n$. Define  the quantity
\begin{equation}
Z_n=\sum_{i=1}^na_i(X_{(i)}-\overline{X}_n). %(3)
\end{equation}
Then the sample mean $\overline{X}_n$ and  the quotient  $Z_n/\overline{X}_{n}$ are
independent iff \ $F$ is a gamma distribution, i.e.,
\[
\overline{X}_n \ \perp\!\!\!\perp \ Z_n/\,\overline{X}_n \quad \Longleftrightarrow \quad F \mbox{ is a gamma distribution}.
\]
} ~~~It is enough to see that $Z_n$ in (3) is a proper feasible
definite statistic of h-degree $1$. If $a_1 = a_2 = \cdots = a_n=a$,
we have $Z_n \equiv 0$, a degenerate case, which is not interesting.

Moreover, when $a_1 =-1,\,  a_2 = \cdots = a_{n-1}=0$ and
$a_n=1$, the quantity $Z_n$ in (3) reduces to the sample range
$R_n.$ Thus Corollary 2 is an extension of Corollary 1.

\smallskip
{\bf Corollary 3.} {\it Given are $n+1$ real numbers,  $p,
a_1, a_2, \ldots, a_{n-1}, a_n$, such that $p>0, a_1>0, a_n>0,$
while $a_{i}\ge 0,$ for $i=2, 3, \ldots, n-1.$  Consider the
quantity
\begin{eqnarray}
Z_n=\sum_{i=1}^na_{i}|X_{(i)}-\overline{X}_n|^p. %(4)
\end{eqnarray}
Then the sample mean $\overline{X}_n$ and  the quotient
$(Z_n)^{1/p}/\overline{X}_n$ are independent iff $F$ is a gamma distribution, i.e.,
\[
\overline{X}_n \ \perp\!\!\!\perp \ (Z_n)^{1/p}/\,\overline{X}_n \quad \Longleftrightarrow \quad F \mbox{ is a gamma distribution}.
\]
} ~~~A little algebra shows that $Z_n$ in (4) is a feasible definite
statistic of h-degree $p.$ Hence  $(Z_n)^{1/p}$ is a feasible
definite statistic of h-order 1, just as required in Theorem 1.

It is worth mentioning that if $p=2$ and $a_i=1/(n-1)$ for all $i,$ the quantity $Z_n$ in (4) reduces
to the sample variance $S_n^2.$ Hence this corollary includes the classical characterization result as a special case.

\medskip
{\bf Corollary 4.}
{\it Suppose that the real numbers $p$ and $a_{i,j}\ge 0,$ where
$1\le i,\ j\le n,$ are such that $p>0$ and $a_{1,n}+a_{n,1}>0.$ Consider the quantity
\begin{eqnarray}
Z_n=\sum_{i=1}^n\sum_{j=1}^na_{i,j}|X_{(i)}-X_{(j)}|^p. %(5)
\end{eqnarray}
Then the sample mean $\overline{X}_n$ and  the quotient
$(Z_n)^{1/p}/\overline{X}_n$ are independent iff $F$ is a gamma distribution, i.e.,
\[
\overline{X}_n \ \perp\!\!\!\perp \ (Z_n)^{1/p}/\,\overline{X}_n \quad \Longleftrightarrow \quad F \mbox{ is a gamma distribution}.
\]
} ~~~We easily check that $Z_n$ is a feasible definite statistic of
h-degree $p,$ and as mentioned above, $(Z_n)^{1/p}$ is such a
statistic of h-degree 1. It remains to refer to Theorem 1.

\medskip
{\bf Corollary 5.} {\it  Consider the well-known {\bf
\emph{Gini's mean difference}}
\begin{equation}
G_n :=\frac{1}{n(n-1)}\sum_{i=1}^n\sum_{j=1}^n|X_i-X_j|. %(6)
\end{equation}
Then the sample mean $\overline{X}_n$ and the quotient
$G_n/\overline{X}_n$ are independent iff
 $F$ is a gamma distribution, i.e.,
\[
\overline{X}_n \ \perp\!\!\!\perp \ G_n/\,\overline{X}_n \quad \Longleftrightarrow \quad F \mbox{ is a gamma distribution}.
\]
} ~~~It is not immediately clear that $G_n$ is a feasible definite
statistic. To see this, we first check that  $G_n$ in (6)  can be
transformed as follows:
\[
G_n
=\frac{1}{n(n-1)}\sum_{i=1}^n\sum_{j=1}^n|X_{(i)}-X_{(j)}|=\frac{4}{n(n-1)}\sum_{i=1}^n\left(i-\frac{n+1}{2}\right)\,X_{(i)}.
\]
For the last representation, see [2], pp. 249, 279. It
 shows that $G_n$ is indeed a feasible definite statistic of h-degree 1, and this is  what we need in  Theorem 1.

Moreover, we easily observe that the quantity $Z_n$ in (5) reduces
to $G_n$ in (6)  when $p=1$ and $a_{i,j}=\frac{1}{n(n-1)},$
 thus Corollary 4 is an extension of Corollary 5. Besides, Corollary 5 also follows from Corollary 2
 due to the last identity for $G_n$ and the fact that $\sum_{i=1}^ni=\frac12n(n+1).$
\medskip\\
\noindent {\bf Corollary 6.} {\it Let $(a_{ij})_{i,j=1}^n$ be a
positive definite $($real$)$ matrix. Consider the quantity
\begin{eqnarray}
Z_n=\sum_{i=1}^n\sum_{j=1}^na_{ij}(X_{(i)}-\overline{X}_n)(X_{(j)}-\overline{X}_n). %(7)
\end{eqnarray}
Then the sample mean $\overline{X}_n$ and the quotient
$(Z_n)^{1/2}/\overline{X}_n$ are independent iff $F$ is a gamma distribution, i.e.,
\[
\overline{X}_n \ \perp\!\!\!\perp \ (Z_n)^{1/2}/\,\overline{X}_n \quad \Longleftrightarrow \quad F \mbox{ is a gamma distribution}.
\]}
 \indent
 Notice that $Z_n$ in (7) is a feasible definite statistic of
h-degree $2,$ hence $(Z_n)^{1/2}$  is a feasible definite statistic
of h-degree 1.

\medskip
{\bf Corollary 7.}
{\it For arbitrary real numbers $p>0, q > 0, a_{ij} \geq 0$ for  $i, j =1, \ldots,n$ with $a_{11} >0$ and $a_{nn}>0, $ consider the quantity
\begin{eqnarray}
Z_n=\sum_{i=1}^n\sum_{j=1}^na_{ij}\,|X_{(i)}-\overline{X}_n|^p\,|X_{(j)}-\overline{X}_n|^q. %(8)
\end{eqnarray}
Then the sample mean $\overline{X}_n$ and  the quotient
$(Z_n)^{1/(p+q)}/\overline{X}_n$ are independent iff $F$ is a gamma distribution, i.e.,
\[
\overline{X}_n \ \perp\!\!\!\perp \ (Z_n)^{1/(p+q)}/\,\overline{X}_n \quad \Longleftrightarrow \quad F \mbox{ is a gamma distribution}.
\]}
 \indent
It follows from (8) that $Z_n$ is a feasible definite statistic of h-degree $p+q,$ which implies that $(Z_n)^{1/(p+q)}$
is a feasible definite statistic of h-degree 1.

\bigskip
{\bf 4. Anosov's theorem, its ramifications and
auxiliary lemmas.}
In order to prove Theorem 1 and see not only the parallel but also
the difference between characterizations of  normal and gamma
distributions, we  recall first Anosov's theorem [1] for normal
distributions (see [8], Chapter 4). We also
provide two ramifications of Anosov's theorem, Theorem A and Theorem
B below, and three auxiliary lemmas.

\smallskip
{\bf Anosov's Theorem.}
{\it Let $n\ge 3$ be an integer and let $X\sim F$ have a
positive continuous density $f$ on $\mathbb R.$ Define the
$(n-2)$-dimensional torus
\begin{eqnarray}
\Phi=\{\phi=(\phi_1,\ldots,\phi_{n-2}):
\phi_j\in[0,\pi],\ j=1,2,\ldots,n-3;\ \phi_{n-2}\in[0,2\pi]\}. %(9)
\end{eqnarray}
Let $T$ be a nonnegative
continuous function on $\Phi$ such that $\int_{\Phi}T(\phi)\,{\rm
d}\phi\in(0,\infty).$ Assume further that $\sigma_j (\phi),\ \phi \in \Phi, \ j=1,\ldots,
n,$ are continuous functions satisfying the following two conditions:
\begin{eqnarray*}
\sum_{j=1}^n\sigma_j(\phi)=0\ \quad  and  \quad
\sum_{j=1}^n\sigma_j^2(\phi)\ \in(0,\infty)\ \ for \ all \ \phi \in \Phi.
\end{eqnarray*}
If  the density $f$ satisfies the integro-functional equation
\begin{eqnarray}
\int_{\Phi}\Pi_{j=1}^n\,f(x+s\sigma_j(\phi))\,T(\phi)\,{\rm d}\phi =
c\,(f(x))^n\!\!
\int_{\Phi}\Pi_{j=1}^nf(s\sigma_j(\phi))\,T(\phi)\,{\rm d}\phi %(10)
\end{eqnarray}
for all $x \in {\mathbb R}$    and all  $s \ge 0$, where $c>0$ is a constant, then $f$ is a normal density function.
}

In order to make one step more, we start with $n-2$ real numbers
$t_1, \ldots, t_{n-2}$ and introduce $n- 2$ functions $g_1, g_2,
\ldots, g_{n-2},$ where  $g_k$ depends on the first $k$ arguments
$t_1, \ldots, t_k,$ as follows:
\[g_{k}(t_1, t_2, \ldots, t_k) =1-t_1^2-t_2^2-\cdots-t_{k}^2, \quad k=1,
2, \ldots, n-2.\]
\noindent The real vector $(t_1, \ldots, t_{n-2})$  belongs to a set
denoted by ${\bf B}_{n-2},$ a subset of ${\mathbb R}^{n-2},$ and we
follow a step-by-step procedure to determine the range of the
variables $t_1, t_2,$ etc., up to $t_{n-2}.$ We start with $t_1$ and
its range $[-1, -\frac{1}{n-1}],$ depending on $n$, hence $g_1(t_1)$
is nonnegative. We continue and find the range for $t_2$, via $n,
t_1$ and $g_1^{1/2}$. Formally we define ${\bf B}_{n-2}$ as follows:
\begin{eqnarray} % (11)
{\bf B}_{n-2}&=&\Big\{(t_1,t_2,\ldots,t_{n-2}):\ -1\le t_1\le -\frac{1}{n-1}, \\ %\nonumber\\
& &\max\Big\{\Big[\frac{n-k+2}{n-k}\Big]^{1/2}\cdot t_{k-1},\
-g_{k-1}^{1/2}\Big\} \le t_k\le \frac{-g^{1/2}_{k-1}}{n-k},\quad k = 2, \ldots, n-2 \Big\}. \nonumber
\end{eqnarray}
%& &~~~\ 2\le k\le n-2\ \Big\}. %(11)
%\end{eqnarray}
It can be checked that the function $g_{n-2}(t),\ t = (t_1, \ldots,
t_{n-2}) \in {\bf B}_{n-2},$ is nonnegative and that the function
${g^{-1/2}_{n-2}}$ is integrable (see Lemma 2 below).

\medskip
In their analysis of the normal distribution, the authors of [3]
applied the following result which is a ramification of the above
Anosov's characterization.

\smallskip
{\bf Theorem A.} {\it If the compact set $\Phi$ in $(9)$ is replaced
by the set ${\bf B}_{n-2}$ in $(11),$ then the conclusion of
Anosov's theorem remains true.}

\smallskip
Let us focus now on tools and ideas allowing us to study gamma
distributions. The first one is an important result from [5], [6];
see Theorem B below. This result is another ramification of Anosov's
theorem, covering a distribution which is not normal. We use the
notation $t:=(t_1, t_2, \ldots, t_{n-2})$ for the points in the
domain ${\bf B}_{n-2}.$

\medskip
{\bf Theorem B.} {\it Let $n\ge 3$ be an integer, and let the random
variable $X$, \ $0<X\sim F,$ have a positive continuous density $f$
on $(0,\infty).$ Let $G$ be a distribution on the set ${\bf
B}_{n-2}$ given in $(11).$  Further, assume that $\sigma_j (t),
j=1,\ldots, n,$ are continuous functions on ${\bf B}_{n-2}$
satisfying two conditions:
\begin{eqnarray*}
\sum_{j=1}^n\sigma_j(t)=0  \quad {\rm and} \quad \sum_{j=1}^n\sigma_j^2(t)\ \in(0,\infty)\quad {for \ all }\
t \in {\bf B}_{n-2} .
\end{eqnarray*}
Suppose that  the density $f$ satisfies the following integro-functional equation:
\begin{eqnarray}
& &\int_{{\bf B}_{n-2}}\Pi_{j=1}^nf(x(s\sigma_j(t)+1))\,{\rm d}G(t)\nonumber\\
&=&{\tilde c}\,(f(x))^n\!\!
\int_{{\bf B}_{n-2}}\Pi_{j=1}^nf(s\sigma_j(t)+1)\,{\rm d}G(t) \quad for \ all \ x> 0 \  and \ all \ s\in [0,s_0),~
\end{eqnarray}
where $\tilde c$ and $s_0$ are two positive constants. Then $f$ is a gamma density.}

\medskip
{\it Remark} 3.
Comparing (10) and (12), we can see the difference between the integro-functional equations in the
normal and the gamma cases. The unknown density $f$ depends on the parameters in quite a different way.
Let us mention that if we further assume in Theorem B that the
density function $f$ is {\it continuously twice-differentiable} on
$(0,\infty),$ then the proof can be essentially simplified by just
differentiating the equation (12) twice (with respect to $s$) at
$s=0,$  and then solving the obtained ordinary differential equation:
\[[f^{\prime\prime}(x)f(x)-(f^{\prime}(x))^2]x^2=c^{\prime}f^2(x),\quad x>0,
\]
where $c^{\prime}\in\R$ is a constant. It is interesting to note
that in such a derivation procedure, the functions $\sigma_j$ and
the integral over the set ${\bf B}_{n-2}$ including the distribution
$G$ in (12) can be canceled out. Thus, they eventually play no roles
in the final equation to be solved. Important is that under the
smoothness and positivity conditions on the density function $f$ on
$(0,\infty),$ the solution $f$ of the above equation is indeed a
gamma density.

 We need also two lemmas,
given below; their proofs are given in [3].

{\bf Lemma 1.} {\it Let $S_n$ be  the sample deviation of a
sample from $F$ of size $n\ge 3$ and let $Z_n$, defined in $(2)$, be
a feasible definite statistic on ${\bf A}$ of h-degree $1.$  Then
there exist two positive constants $c_*<c^*$ such that
\[
0<c_* \le \frac{Z_n}{S_n} \le c^*<\infty\ \  a.s.
\]
}
{\bf Lemma 2.} {\it Consider the function
$g_{n-2}(t)=1-t_1^2-t_2^2-\cdots-t_{n-2}^2$ introduced above.  Then
\[
I_{n-2} := \int_{{\bf B}_{n-2}}\frac{{\rm d}t_1 \cdots {\rm d}t_{n-2}}{\sqrt {g_{n-2}(t)}}
=\frac{2\,\pi^{(n-1)/2}}{n!\Gamma((n-1)/2)}.
\]}

Since $g_{n-2}$ is nonnegative and $I_{n-2}< \infty$, the normalized  function
\[
\tilde{g}(t)=\frac{1} {I_{n-2}\,\sqrt{g_{n-2}(t)}}, \quad t \in  {\bf B}_{n-2},
\]
is a proper density function on the set ${\bf B}_{n-2}$ (see [4]). Similarly,
Lemmas 1 and 2 together will be used to define another distribution
$G$ required in Theorem B when proving Theorem 1 (see relation (26) below).

As in the beginning of Section 3,  let $X_{(1)}\le X_{(2)}\le
\cdots\le X_{(n)}$ be the order statistics of a random sample $X_1,
X_2,\ldots, X_n$ of size $n\ge 3$ from a distribution $F$ which has
positive continuous density $f$ on $(0,\infty).$ The corresponding
realized ordered values are   $x_{(1)}\le x_{(2)}\le \cdots\le
x_{(n)},$ with their sample mean and  sample variance:
\[
\overline{x}_n=\frac1n\sum_{i=1}^nx_{(i)},\quad
s_n^2=\frac{1}{n-1}\sum_{i=1}^n(x_{(i)}-\overline{x}_n)^2.
\]

Following [4], we define the transformation
\begin{eqnarray*}T^*_1:(x_{(1)}, x_{(2)}, \ldots, x_{(n)}) \ \mapsto \ (t_1,t_2,\ldots,t_{n-2},w_1,w_2)\end{eqnarray*} as follows:
\begin{eqnarray*}
\begin{cases}\displaystyle
t_i=\Big[\frac{n-i+1}{(n-1)(n-i)}\Big]^{1/2}\Big[\frac{x_{(i)}-\overline{x}_n}{s_n}+\frac{1}{n-i+1}\sum_{k=1}^{i-1}\frac{x_{(k)}-\overline{x}_n}{s_n}\Big],\ \
1\le i\le n-2, \vspace{0.2cm}\\
w_1=\overline{x}_n, \vspace{0.1cm}\\
 w_2=s_n.\\
\end{cases}
\end{eqnarray*}
The summation above in $t_i$ is taken to be zero if $i=1.$ Then we
use the functions $g_k(t_1, \ldots, t_k)=1-\sum_{i=1}^k t_i^2, \ 1
\le k\le n-2,$ described before, and denote by $(T_1^*)^{-1}$ the
inverse transformation to $T_1^*$:
\begin{eqnarray*}(T_1^*)^{-1}:
(t_1,t_2,\ldots,t_{n-2},w_1,w_2) \ \mapsto \  (x_{(1)}, x_{(2)},
\ldots, x_{(n)}).
\end{eqnarray*}
$(T_1^*)^{-1}$ is defined through the following relations:
\begin{eqnarray}
\begin{cases}\displaystyle
\frac{x_{(i)}-w_1}{w_2\sqrt{n-1}}=\Big[\frac{n-i}{n-i+1}\Big]^{1/2}\cdot t_i
-\sum_{k=1}^{i-1}\frac{t_k}{[(n-k)(n-k+1)]^{1/2}},\ \ 1\le i\le n-2, \vspace{0.2cm}\\\displaystyle
\frac{x_{(n-1)}-w_1}{w_2\sqrt{n-1}}=-\sum_{k=1}^{n-2}\frac{t_k}{[(n-k)(n-k+1)]^{1/2}}-[g_{n-2}/2]^{1/2}, \vspace{0.1cm}\\\displaystyle
\frac{x_{(n)}-w_1}{w_2\sqrt{n-1}}=-\sum_{k=1}^{n-2}\frac{t_k}{[(n-k)(n-k+1)]^{1/2}}+[g_{n-2}/2]^{1/2},\\
\end{cases}
\end{eqnarray}
where \ $\overline{x}_n=w_1$ and \ $s_n=w_2.$

Further, consider also another transformation, namely,
\begin{eqnarray}
T_2^*:(x_{(1)}, x_{(2)}, \ldots,
x_{(n)}) \ \mapsto \  (t_1,t_2,\ldots,t_{n-2},w_1,w_3),
\end{eqnarray}
where $w_3=w_2/w_1,$ and for $n\ge 3,$ define two more  subsets
${\bf D}^*_n$ and ${\bf R}^*_n$ of ${\mathbb R}^n$:
\begin{eqnarray}
{\bf D}^*_n&=&\{(x_{(1)}, x_{(2)},\ldots, x_{(n)}):\ 0<x_{(1)}\le  x_{(2)}\le \cdots\le  x_{(n)}\};\\
{\bf R}^*_n&=&\Big\{(t_1, t_2,\ldots, t_{n-2}, w_1, w_3):\ \max\Big\{-1, \ -\frac{\sqrt{n}}{(n-1)w_3}\Big\}\le t_1 \le
-\frac{1}{n-1},\nonumber\\
& &~~~\max\Big\{\Big[\frac{n-k+2}{n-k}\Big]^{1/2}\cdot t_{k-1}, \
-g_{k-1}^{1/2}\Big\} \le t_k \le -\frac{g^{1/2}_{k-1}}{n-k},\nonumber\\
& &~~~\ 2\le k \le n-2,\ \ w_1>0,\ \ 0<w_3 \le \sqrt{n}\Big\}.
\end{eqnarray}

We need one result more, Lemma 3 below, which is Theorem 2.2 in
Hwang and Hu (1994), rewritten here in terms of the above notations.

{\bf Lemma 3.} {\it  For $n\ge 3,$ let $T_2^*$ be the
transformation  defined in $(14)$. Then $T_2^*$ establishes a
one-to-one correspondence between  the domain ${\bf D}^*_n$ and the
range ${\bf R}^*_n,$ defined in $(15)$ and $(16)$, respectively,
except for a set of $n$-dimensional Lebesgue measure zero.
Furthermore, the absolute value of the Jacobian of this
transformation is
\begin{eqnarray*}|J_2^*|= \sqrt{n}\,(n-1)^{(n-1)/2}\cdot
w_1^{n-1}w_3^{n-2}\cdot g_{n-2}^{-1/2}.
\end{eqnarray*}}

\smallskip
{\bf 5. Proof of Theorem 1.}

{\it Proof of Theorem 1.} {\it Step 1. Sufficiency.} Since
the degree of homogeneity of the feasible definite statistic  $Z_n$
is 1, it follows that:
\begin{eqnarray*}
Z_n&=&U(X_{(1)}-\overline{X}_n,X_{(2)}-\overline{X}_n,\ldots,X_{(n)}-\overline{X}_n)\nonumber\\
&=&\overline{X}_nU((X_{(1)}-\overline{X}_n)/\overline{X}_n,(X_{(2)}-\overline{X}_n)/\overline{X}_n,\ldots,(X_{(n)}-\overline{X}_n)/\overline{X}_n)\nonumber\\
&=&\overline{X}_nU(X_{(1)}/\overline{X}_n-1,X_{(2)}/\overline{X}_n-1,\ldots,X_{(n)}/\overline{X}_n-1).
\end{eqnarray*}
Hence,
\begin{eqnarray}
\frac{Z_n}{\overline{X}_n}=U\left(\frac{X_{(1)}}{\overline{X}_n}-1, \frac{X_{(2)}}{\overline{X}_n}-1, \ldots, \frac{X_{(n)}}{
\overline{X}_n}-1\right).
\end{eqnarray}
 The RHS of relation (17) is a function of the random vector
 $\left(\frac{X_{(1)}}{\overline{X}_n}, \frac{X_{(2)}}{\overline{X}_n}, \ldots, \frac{X_{(n)}}{\overline{X}_n}\right).$
If $F$  is a gamma distribution, then $\overline{X}_n$ is
independent of $\left(\frac{X_{(1)}}{\overline{X}_n},
\frac{X_{(2)}}{\overline{X}_n}, \ldots,
\frac{X_{(n)}}{\overline{X}_n}\right)$ and hence $\overline{X}_n$ is
independent of ${Z_n}/{\overline{X}_n}.$ See also [10] or [9].

This proves the sufficiency part of the theorem.

{\it Step 2. Necessity.}  Suppose that the sample mean
$\overline{X}_n$ of a sample of size $n$ coming from a distribution
$F$ is independent of the quotient $V_n=Z_n/\overline{X}_n,$ where
$Z_n$ is a homogeneous feasible definite   statistic of degree 1.
Then we want to show that $F$ is a gamma distribution. By the
homogeneity property again, rewrite
\begin{eqnarray*}
Z_n&=&U(X_{(1)}-\overline{X}_n,X_{(2)}-\overline{X}_n,\ldots,X_{(n)}-\overline{X}_n)
= S_n\,U(\Lambda_1,\Lambda_2,\ldots,\Lambda_{n}),
\end{eqnarray*}
where
\begin{eqnarray}
\Lambda_i=\frac{X_{(i)}-\overline{X}_n}{S_n},\quad
\ i=1, 2, \ldots, n,
\end{eqnarray}
and $(\Lambda_1,\Lambda_2,\ldots,\Lambda_{n})$ takes a.s. values
$(\tilde\lambda_1,\tilde\lambda_2,\ldots,\tilde\lambda_{n})$ in the
compact subset ${\bf A}_n\subset {\bf A}$,
\begin{eqnarray*}
{\bf A}_n=\Big\{\lambda=(\lambda_1,\ldots,\lambda_n): \lambda_1\le \lambda_2\le\cdots\le \lambda_n, \
\small\sum_{i=1}^n\lambda_i=0,\ \sum_{i=1}^n\lambda_i^2=n-1\Big\}.
\end{eqnarray*}
 It then follows from Lemma 1 that there exist two
positive constants $c_*<c^*$ such that
\begin{eqnarray}
0<c_*\le \frac{Z_n}{S_n}=U(\Lambda_1,\Lambda_2,\ldots,\Lambda_{n})\le c^*<\infty \   \mbox{ a.s.}
\end{eqnarray}

Write the realized values of
\[S_n=\frac{Z_n}{U(\Lambda_1,\Lambda_2,\ldots,\Lambda_{n})}
\]
as
\begin{eqnarray}
s_n=\frac{z_n}{U(\tilde\lambda_1,\tilde\lambda_2,\ldots,\tilde\lambda_{n})}:=\frac{z_n}{U(\tilde\lambda)},
\quad \tilde\lambda=(\tilde\lambda_1,\ldots,\tilde\lambda_n),
\end{eqnarray}
where each $\tilde\lambda_i$ is a function of
$t=(t_1,t_2,\ldots,t_{n-2})$ defined through (13) and (18).

Recall now that the joint density of the order statistics $(X_{(1)},
X_{(2)}, \ldots, X_{(n)}),$ denoted by $h$, is expressed via the
1-dimensional density $f$ of $X$:
\[
h(x_{(1)}, x_{(2)}, \ldots,
x_{(n)})=n!\prod_{i=1}^nf(x_{(i)}),\ \ 0<x_{(1)}\le x_{(2)}\le \cdots\le
x_{(n)}.
\]
Next, consider the composition of two transformations, $T_2^*$ from
(14) and a new one, $T^*,$ defined by
\[
T^*: \ (t_1,t_2,\ldots,t_{n-2},w_1,w_3) \ \mapsto \ (t_1,t_2,\ldots,t_{n-2},w_1,v_n).
\]
Here, $w_1=\overline{x}_n,\  w_3=s_n/w_1,$
$v_n=z_n/w_1=U(\tilde\lambda)w_3$ and $z_n$ is defined in (20). The
transformation $T^*$ has a Jacobian  $J^*=1/U(\tilde\lambda)$ which
is positive. Hence, the composition $T_2^*\circ T^*$ has a Jacobian
with absolute value (by Lemma 3):
\begin{eqnarray*}
|J_2^*|J^*&=& \sqrt{n}\,(n-1)^{(n-1)/2}\cdot
w_1^{n-1}w_3^{n-2}\cdot g_{n-2}^{-1/2}(t)\cdot \frac{1}{U(\tilde\lambda)}\\
&=&\sqrt{n}\,(n-1)^{(n-1)/2}\cdot(\overline{x}_n)^{n-1}
v_{n}^{n-2}\Big(\frac{1}{U(\tilde\lambda)}\Big)^{n-1}\cdot g_{n-2}^{-1/2}(t),
\end{eqnarray*}
in which, to recall that $g_{n-2}(t)=1-\sum_{i=1}^{n-2}t_i^2.$ Therefore, the joint
density of the statistics $(T_1,T_2,\ldots,T_{n-2},$
$\overline{X}_n,V_n)$ is of the form:
\[
\tilde{h}(t_1,t_2,\ldots,t_{n-2},\overline{x}_n,v_n)=n!\prod_{i=1}^nf(x_{(i)})|J_2^*|J^*,
\]
where each $T_i$, a function of order statistics $X_{(1)},
\ldots,X_{(n)}$ (or of statistics $\Lambda_1,\ldots,\Lambda_{n-2}$
in (18)), takes values $t_i$ defined in the transformation $T_1^*.$

Notice now that by (13), (18)  and (20), we have, for $i=1, 2, \ldots, n,$ that
\[
x_{(i)}=\overline{x}_n+\frac{v_n\overline{x}_n}{U(\tilde\lambda)}\cdot\tilde\lambda_i(t_1,t_2,\ldots,t_{n-2})
=\overline{x}_n\Big[1+\frac{v_n}{U(\tilde\lambda)}\cdot\tilde\lambda_i(t_1,t_2,\ldots,t_{n-2})\Big].\quad
\]

Denote $v_0={c_*\sqrt{n}}/{(n-1)}.$  If $v_n\in(0,v_0),$ then we
have, by (19) and (20), that
\[
w_3 \le \frac{c_*}{U(\tilde\lambda)}\frac{\sqrt{n}}{n-1} \le \frac{\sqrt{n}}{n-1},
\]
and hence $(t_1,t_2,\ldots,t_{n-2})\in {\bf B}_{n-2}$,  provided
that $(t_1,t_2,\ldots,t_{n-2}, w_1,w_3)\in{\bf R}^*_{n}.$ This helps
us to simplify the domain of the joint density
$\tilde{h}(t_1,t_2,\ldots,t_{n-2},\overline{x}_n,v_n)$ when
$v_n\in(0,v_0).$ Thus, for $\overline{x}_n>0$ and $v_n\in(0,v_0),$
the 2-dimensional marginal joint density of $(\overline{X}_n,V_n)$
has the following form (for simplicity, we use the notations
$\overline{x}=\overline{x}_n, \ v=v_n$):
\begin{eqnarray} %(21)
\tilde{h}(\overline{x},v)&=&n!\sqrt{n}\,(n-1)^{(n-1)/2}\cdot(\overline{x})^{n-1}v^{n-2}\nonumber \\
& & \times \int_{{\bf B}_{n-2}}\Big(\frac{1}{U(\tilde\lambda)}\Big)^{n-1}
\, g_{n-2}^{-1/2}(t)\prod_{i=1}^n
f\Big(\overline{x}\big[1+\frac{v}{U(\tilde\lambda)}\,\tilde\lambda_i(t)\big]\Big)\,{\rm d}t_1 \cdots {\rm d}t_{n-2},~~
\end{eqnarray}
where $t=(t_1,t_2,\ldots,t_{n-2})\in{\bf B}_{n-2}\subset{\mathbb
R}^{n-2}$ defined in (11).

Now we apply the independence condition on $\overline{X}_n$ and
$V_n,$ and write the joint density in (21) as the product of the
densities of $\overline{X}_n$ and $V_n:$
\begin{eqnarray}
\tilde{h}(\overline{x},v)=\tilde{h}_{\overline{X}_n}(\overline{x})\,\tilde{h}_{V_n}(v), \quad {\rm if} \quad \overline{x}>0 \quad {\rm and}
\quad v\in(0,v_0).
\end{eqnarray}
Letting $\overline{x}=1$ in (21) and (22), we obtain the density of
$V_n:$
\begin{eqnarray}
\tilde{h}_{V_n}(v)&=&\frac{1}{\tilde{h}_{\overline{X}_n}(1)}\,\tilde{h}(1,v)\nonumber\\
&=& \frac{1}{\tilde{h}_{\overline{X}_n}(1)}\,n!\sqrt{n}\,(n-1)^{(n-1)/2}\cdot
v^{n-2}\int_{{\bf B}_{n-2}}\Big(\frac{1}{U(\tilde\lambda)}\Big)^{n-1}\nonumber\\
&  &\times\, g_{n-2}^{-1/2}(t)\prod_{i=1}^n
f\Big(1+\frac{v}{U(\tilde\lambda)}\cdot\tilde\lambda_i(t)\Big)\,{\rm d}t_1 \cdots {\rm d}t_{n-2}, \quad {\rm if} \   v\in(0,v_0).
\end{eqnarray}
Plugging (21) and (23) in (22), canceling the common term
$n!\sqrt{n}\,(n-1)^{(n-1)/2}\cdot v^{n-2}$, then letting $v\to
0$ and  canceling the common integral term, we find the density of
$\overline{X}_n:$
\begin{eqnarray}
\tilde{h}_{\overline{X}_n}(\overline{x})={\tilde c}\,(\overline{x})^{n-1}\,(f(\overline{x}))^n, \quad \overline{x}>0; \quad
\tilde{h}_{\overline{X}_n}(\overline{x}) = 0, \ \mbox{ otherwise,}
\end{eqnarray}
where ${\tilde c}=\tilde{h}_{\overline{X}_n}(1)/(f(1))^n > 0$ is a constant.

Combining (21) through (24), we finally have the integro-functional
equation\,:
\begin{eqnarray}
& &\int_{{\bf B}_{n-2}}\Big(\frac{1}{U(\tilde\lambda)}\Big)^{n-1}\cdot g_{n-2}^{-1/2}(t)
\prod_{i=1}^n
f\Big(\overline{x}\big(1+\frac{v}{U(\tilde\lambda)}\cdot\tilde\lambda_i(t)\big)\Big)\,{\rm d}t_1 \cdots {\rm d}t_{n-2}\nonumber\\
&=&\!\!C\,(f(\overline{x}))^n\!\!\int_{{\bf B}_{n-2}}\!\Big(\frac{1}{U(\tilde\lambda)}\Big)^{n-1}\!\cdot g_{n-2}^{-1/2}(t)
\prod_{i=1}^n
f\Big(1+\frac{v}{U(\tilde\lambda)}\cdot\tilde\lambda_i(t)\Big)\,{\rm d}t_1 \cdots {\rm d}t_{n-2},~~~
\end{eqnarray}
for all $\overline{x}>0$ and all $v\in[0,v_0),$ where
$C=1/(f(1))^n>0$ is a constant.

It is seen that (25) is of the form (12) in Theorem B,  because, by
Lemmas 1 and 2 and Eq.\,(18), we have
\[
M:=\int_{{\bf B}_{n-2}}\Big(\frac{1}{U(\tilde\lambda)}\Big)^{n-1}\cdot g_{n-2}^{-1/2}(t)\,{\rm d}t_1 \cdots {\rm d}t_{n-2} \ \in \ (0,\infty),
\]
\[\sum_{i=1}^n\sigma_i(t):=\sum_{i=1}^{n}\frac{\tilde\lambda_i(t)}{U(\tilde\lambda)}=0\ \ \mbox{ for all } \ t\in {\bf B}_{n-2},
\]
\[\sum_{i=1}^n\sigma_i^2(t)=\sum_{i=1}^{n}\Big[\frac{\tilde\lambda_i(t)}{U(\tilde\lambda)}\Big]^2
=\frac{n-1}{U^2(\tilde\lambda)}\ \ \in  (0,\infty) \ \ \mbox{ for all } \ t\in {\bf B}_{n-2}.
\]

Here, the last two conditions are required in Theorem B, and we may
define the distribution $G(t)$ for  $t=(t_1,t_2,\ldots,t_{n-2}) \in {\bf B}_{n-2},$ as follows:
\begin{equation}
{\rm d}G(t) := \frac{1}{M}\Big(\frac{1}{U(\tilde\lambda)}\Big)^{n-1}\cdot
g_{n-2}^{-1/2}(t)\, {\rm d}t_1\, \cdots {\rm d}t_{n-2}.
\end{equation}

Notice that here, $\overline{x}$ and $v$ play the roles of $x$ and
$s$ in Eq.\,(12), respectively. Therefore, by Theorem B, $f$ is a gamma density. This completes the proof of Theorem 1.\ \ $\qed$

\bigskip
{\bf 6. Final comments and two conjectures.}
In our main results, Theorem 1 and Corollaries 1 -- 7, we assume that the
sample size is $n\ge 3$ and impose smoothness condition on the
underlying distribution $F$ of $X.$ A natural question arises:

\smallskip
{\bf Question.} \ {\it Is it
possible to relax, or dispose, the smoothness condition on $F\,?$ }

Here are some motivating arguments when looking closely to the case of a sample size $n=2.$
In this case we have a simple relationship between the sample variance $S_2^2$ and the sample range $R_2,$
namely $S_2^2=\frac12R_{2}^2.$ Hence, if the sample mean $\overline{X}_2$
and the quotient $R_{2}/\overline{X}_2$ are independent, then so are
$\overline{X}_2$ and $S_2^2/\overline{X}_2^2.$ This in turn implies
that $\overline{X}_2$ and $S_2/\overline{X}_2$ are independent. The latter property,
 according to  the classical result, characterizes the
underlying distribution $F$ as being a gamma  distribution. Notice
that this conclusion holds without any smoothness condition on $F$.
Similarly, for Gini's mean difference in (6), we find that $G_2$
exactly equals $R_2.$ Therefore, if $\overline{X}_2$ and
$G_2/\overline{X}_2$ are independent, then $F$ is a gamma
distribution.

These observations lead to the idea of eventually extending some of our results, in
particular, the statements in Corollary 1 and Corollary 5, without assuming any smoothness condition
on the underlying distribution $F.$

{\bf Conjecture 1.} {\it Let $X_1, X_2,\ldots, X_n$ be a
random sample of size $n\ge 3$ from {\bf \emph{any}} nondegenerate
distribution $F$ on $(0,\infty).$  If the sample mean
$\overline{X}_n$ and the quotient $R_{n}/\overline{X}_n$ are
independent, then $F$ is a gamma distribution. }

{\bf Conjecture 2.} {\it Let $X_1, X_2,\ldots, X_n$ be a
random sample of size $n\ge 3$ from {\bf \emph{any}} nondegenerate
distribution $F$ on $(0,\infty).$  If the sample mean
$\overline{X}_n$ and the quotient $G_n/\overline{X}_n$  are
independent, then $F$ is a gamma distribution.}

{\bf Acknowledgment.} \ We are grateful to Prof. Chin-Yuan
Hu (National Changhua University of Education, Taiwan) for the
fruitful discussions and collaboration over the years, and in
particular for his useful comments when preparing this paper.

\bigskip
\begin{center}
REFERENCES
\end{center}

\noindent
1. D.\,V. Anosov,
``On an integral equation from statistics", {\it Vestnik Leningrad State Univ., Ser. III Math. Mech.
Astron.} \ {\bf 19} (1964), 151--154. (In Russian)

\noindent
2. H.\,A. David, H.\,N. Nagaraja,  {\it Order Statistics,} 3rd ed., Wiley, Hoboken (NJ), 2003.

\noindent
3. C.\,Y. Hu, G.\,D. Lin,  ``Characterizations of the
normal distribution via the independence of the sample mean and
the feasible definite statistics with ordered arguments", {\it
Ann. Inst. Statist. Math.}  {\bf 74} (2022), 473--488.

\noindent
4. T.-Y. Hwang, C.-Y. Hu,  ``On the joint
distribution of Studentized order statistics",  {\it Ann. Inst. Statist. Math.} \  {\bf 46} (1994), 165--177.

\noindent 5.  T.-Y. Hwang, C.-Y. Hu,  ``On an integro-functional
equation arising in order statistics", {\it Bull. Inst. Math. Acad.
Sinica} \ {\bf 25} (1997), 161--169.

\noindent
6.  T.-Y. Hwang, C.-Y. Hu,  ``On a characterization of the gamma distribution: the independence
of the sample mean and the sample coefficient of variation",
{\it Ann. Inst. Statist. Math.} \ {\bf 51} (1999), 749--753.

\noindent
7.  T.-Y. Hwang, C.-Y. Hu,  ``On some characterizations of population distributions",
{\it  Taiwanese J. Math.} \  {\bf 4} (2000),  427--437.

\noindent
8. A.\,M.\, Kagan, Yu.\,V.  Linnik, C.\,R. Rao,  {\it Characterization Problems in
 Mathematical Statistics.}  Wiley, New York, 1973.

\noindent
9.  G.\, Marsaglia,  ``Extension and applications of Lukacs' characterization of the gamma distribution", in:
{\it Proc. Symp. Statist. Related Topics}, Ottawa, Oct. 24--26, 1974, Ed. A.K.M.E. Saleh,  9.01--9.13, Carleton University, Ottawa.

\noindent
10. E.\,J.\,G. Pitman,    ``The `closest estimates' of statistical parameters",
{\it Math. Proc. Cambr. Philos. Soc.}  {\bf 44} (1937), 212--222.

\end{document}